\documentclass[12pt,twoside]{amsart}
\usepackage[margin=1.25in]{geometry}

\usepackage{mathpazo}
\usepackage{tgbonum}
\usepackage{comment}

\usepackage[T1]{fontenc}

\usepackage{physics}
\usepackage{amsmath}
\usepackage{amsthm}
\usepackage{amsfonts}
\usepackage{amssymb}
\usepackage{esvect}
\usepackage[dvipsnames]{xcolor}
\usepackage{graphicx}
\usepackage{float}
\usepackage[margin=2pc,font=small,labelfont=bf]{caption}
\usepackage{relsize}
\usepackage{bigints}
\usepackage[makeroom]{cancel}
\usepackage{ wasysym }
\usepackage{faktor}
\usepackage{bbm}
\usepackage{cancel}
\usepackage{tikz-cd}
\usepackage{csquotes}
\usepackage[english]{babel}
\usepackage[hidelinks]{hyperref}
\usepackage{mathtools}
\usepackage{multicol}
\usepackage{biblatex}%bib
\usepackage{wrapfig}
\usepackage{cleveref}

\usepackage{rotating}
\usepackage{shortcuts}
\usepackage{braket}

\usepackage{todonotes}

\DeclareFieldFormat{pages}{#1}
\DeclareFieldFormat[article]{pages}{#1}

\graphicspath{ {./img/} }

\usetikzlibrary{hobby}

\title{
{Inducing coverings on Hilbert schemes}
}

\author{Lucas Li Bassi}
\address{Lucas Li Bassi, 	DIMA (Università di Genova), Via Dodecaneso 35, 16146 Genova, Italy}
\email{lucas.libassi@gmail.com}
\urladdr{http://https://sites.google.com/view/lucaslibassi/}

\author{Filippo Papallo}
\address{Filippo Papallo, 	DIMA (Università di Genova), Via Dodecaneso 35, 16146 Genova, Italy}
\email{papallo@dima.unige.it}
\urladdr{}

\date{}

\begin{document}

\begin{abstract}
    We find an explicit geometric description of all coverings of $\hilb\Sigma$ when $\Sigma$ is a normal, complex, quasi-projective surface with finite fundamental group. We then apply this construction to show that if $\Sigma$ is an irreducible symplectic surface then $\hilb\Sigma$ is an irreducible symplectic variety.
\end{abstract}
\maketitle

\tableofcontents

    \section*{Introduction} 
        
        Since their discovery, irreducible holomorphic symplectic manifolds (or ISM) have attracted the interest of the mathematical community. These varieties are compact % complex 
        Kähler manifold X that are simply connected and such that $H^0(X, \Omega^2_X)$ is generated by an everywhere non-degenerate holomorphic 2-form, which is referred to as the \emph{symplectic} form. The first examples of such manifolds in dimension higher than two were constructed by \textbf{Beauville} in his seminal paper \parencite[Deuxième~partie]{BEAUVILLE1983} using the construction of Hilbert schemes of $n$ points
        $\operatorname{Hilb}^{n}(\Sigma)$ on a surface $\Sigma$
        (this surface being an abelian or a $K3$ surface).  
        Their importance is further underscored by the very well known \textbf{BB decomposition}.
        \begin{thm*}[\textbf{Beauville--Bogomolov decomposition}]
            Let $X$ be a smooth connected compact K\"ahler manifold with $c_1(X) = 0$. There exists a finite étale cover of X which is a product of a complex tori, irreducible symplectic manifolds and irreducible Calabi-Yau manifolds.
        \end{thm*}
        In recent years, attention has shifted towards extending this theory to the singular settings. There are many ways to generalize the notion of ISM in order to allow singularities. Three of the most used are surely the notions of \emph{irreducible symplectic orbifold} (ISO), \emph{irreducible symplectic variety} (ISV) and \emph{primitive symplectic variety} (PSV) and they are related by the chain of implications $$\text{\textbf{ISO}} \implies 
            \text{\textbf{ISV}} \implies \text{\textbf{PSV}}.$$ We refer to \cite[Section 1]{GPP-SSS} for definitions and a detailed overview. In that article, \textbf{Garbagnati--Penegini--Perego} \parencite[Theorem~{1.8}]{GPP-SSS} classified all projective symplectic surfaces. Moreover, they were able to prove that, if the surface $\Sigma$ is a PSV (resp. ISO), then $\hilb{\Sigma}$ is a PSV (resp. ISO). However, they leave open the case where $\Sigma$ is an ISV. The goal of this paper is to bridge this gap. We prove:
       \begin{thmx}\label{ThmB}
            Let $\Sigma$ be a projective irreducible symplectic surface,
            then $\hilb{\Sigma}$ is an irreducible symplectic variety
            of dimension $4$.
        \end{thmx}       
    
In order to verify the definition of ISV (\Cref{ISV}) a crucial point is to be able to control, for every finite quasi-\'etale covering, the behavior of the pull-back of the symplectic form. This led to the more general question of explicitly constructing all (quasi-)étale coverings of the Hilbert scheme of two points on a normal surface. Indeed, at the moment we are writing this article, there exist no general method to induce from a (quasi-)étale covering of a surface $\Sigma$ a (quasi-)étale covering of $\hilb\Sigma$. Such a result can be found in \cite{TARO} when $\Sigma$ is an Enriques surface. We were able to modify his construction to answer our question and we proved the following theorem.
        \begin{thmx} \label{ThmA}
            Let $\Sigma$ be a normal quasi-projective complex surface with ADE singularities and finite $\pi_1(\Sigma)$. 
            Then, every (quasi-)étale cover of $\hilb\Sigma$ is induced by a (quasi-)étale cover of $\Sigma$.
        \end{thmx}

 \paragraph{Plan of the paper.} In \Cref{preliminaries} we recall the necessary background on (quasi-)étale maps and on the geometry of Hilbert schemes of points on normal surfaces. In \Cref{cover} we construct all (quasi-)étale coverings of the Hilbert scheme of two points on a normal quasi-projective surface with finite fundamental group, showing that they are induced by coverings of the surface itself. \Cref{sec:ISV} is devoted to the study of symplectic structures: we recall the notion of irreducible symplectic variety and apply the results of \Cref{cover} to prove that if a surface is irreducible symplectic, then so is its Hilbert square.

    \section*{Acknowledgments}

        We wish to thank 
        A. Perego for suggesting us such a nice
        problem and for his helpful insights on the subject; 
        we would also like to acknowledge V. Bertini, E. Fatighenti and A. Rapagnetta for the helpful discussions. A warm appreciation goes to A. Garbagnati for carefully reading the preliminary version of this work
        and pointing out many mistakes therein. Moreover, we are deeply grateful to A. Frassineti for his continuous 
        encouragement and fruitful conversations 
        during the preparation of this paper.
        The project was supported by the Research Project PRIN 2020 - CuRVI, CUP J37G21000000001. Our gratitude goes the math department of Università di Genova for the fundings and for the nice environment offered during the research.

    \section{Preliminaries}\label{preliminaries}
        
        In this section, we recall the main definitions of coverings 
        we are interested in, both in the topological and in the algebraic setting;
        our conventions follow those by \textbf{Greb-Kebekus-Peternell} \parencite[]{GKP2013}. 
        Then, we review the construction of the Hilbert scheme.
        Throughout this paper, we will always
        deal with normal, connected, quasi-projective %K\"ahler 
        varieties over $\C$;
        these hypothesis allow us to apply Zariski's Main Theorem,
        a very important ingredient for our constructions.

        \subsection{Quasi-\'etale maps and branched coverings}

        The main aim of the paper is to build covering spaces
        of some complex varieties of dimension $4$; 
        there are many different notions
        of covering, according to the category we are working in:
        either solely topological spaces or 
        algebraic varieties.
        Working over the complex numbers allows us to 
        switch among the Zariski and the analytic topology, 
        thus we recall some different definitions of covering.
    
        \begin{df}[Topological covering]
            Given $X$ a topological space, a \textbf{topological covering space} is
            a topological space $Z$, endowed with a
            continuous surjective map $\gamma : Z \to X$
            with the property that, for each $x \in X$, 
            there exists an open neighborhood $U_{x} \subset X$ 
            such that $\gamma^{-1}(U_{x}) = \sqcup_{i \in I} V_{i}$
            is a disjoint union of open subsets that are mapped
            homeomorphically to $U_{x}$, i.e. for each $i \in I$,
            \begin{equation*}
                \gamma\vert_{V_{i}}: V_{i} \overset{\sim}{\longrightarrow} U_{x}
            \end{equation*}
            is a homeomorphism.
            The covering space is called \textbf{finite} if the set of indices $I$ 
            above is finite for every $x \in X$.
        \end{df}

        \begin{df}[Algebraic cover]
            Given a normal, connected variety $X$,
            by \textbf{finite algebraic cover} we mean a finite, surjective
            morphism $\xi : Z \to X$, where $Z$ is a normal,
            connected variety as well.
        \end{df}

        In both algebraic and topological case, given two covers
        $f:Z \to X$ and $f':Z'\to X$, a \textbf{cover isomorphism}
        is an isomorphism (in the appropriate setting) $\phi : Z \overset{\sim}{\to} Z'$
        such that $f = f'\circ \phi$:
        \begin{equation*}
            \begin{tikzcd}
                Z \ar[rr, "\phi"] \ar[dr, "f"'] & & Z' \ar[dl, "f'"] \\
                & X & \,.
            \end{tikzcd}
        \end{equation*}

        \begin{rmk}
            If $\xi : Z \to X$ is a finite algebraic cover, then it defines
            a \emph{connected} topological finite cover if and only if 
            $\xi$ is \textbf{\'etale}, i.e. smooth and unramified.
        \end{rmk}

        When dealing with coverings of \emph{singular} varieties, 
        asking for the cover map to be
        \'etale may be too restrictive. Thus,
        it comes more natural to allow some small enough %(\textcolor{red}{Qui perché emph? Controllare dato che small significa già qualcosa di preciso per mappe})
        ramification: the following definition is due to Catanese 
        \parencite[Definition~{1.1}]{CATANESE-QED}.
        
        \begin{df}[Quasi-\'etale morphism]
            A map of schemes $\epsilon : Z \to X$ is \textbf{quasi-\'etale} if it is 
            quasi-finite and \emph{\'etale in codimension $1$},
            i.e. there exists a closed subset $R \subset Z$ of $\textrm{codim}_{Z}(R) \ge 2$
            such that $\epsilon\vert_{Z \setminus R} \, : \, Z \setminus R \to X$
            is \'etale.
        \end{df}

        We will mainly need the following properties of quasi-\'etale maps:
        
        \begin{prop}
            Let $\xi : Z \to X$ be a \emph{quasi-\'etale} finite (algebraic) cover.
            \begin{enumerate}
                \item If the base $X$ is smooth, then $\xi$ is in fact \'etale,
                and hence a finite topological cover.
                \item If the base $X$ has at worst canonical singularities,
                so does its cover $Z$.
            \end{enumerate}
            \begin{proof}
                For details, we refer to \parencite[Section~{2}]{CATANESE-QED}.
            \end{proof}
        \end{prop}

        Recall that, for any reasonably good enough\footnote{That is, path-connected, locally path-connected and semi-locally simply connected,
        as in the hypothesis of \parencite[Theorem~{1.38}]{HATCHER}. Complex varieties of this paper, seen with the analytic topology, always have these properties.} topological space $X$, its fundamental group $\pi_{1}(X)$ classifies
        all of its coverings: more precisely, 
        to any conjugation class of a subgroup $H \le \pi_{1}(X)$,
        there exists a path-connected topological cover $\xi : Z \to X$ 
        s.t. $\xi_{*}\pi_{1}(Z) = H$, and conversely any such cover
        defines a unique subgroup $\xi_{*}\pi_{1}(Z) \le \pi_{1}(X)$ up to conjugation. Moreover, there is a bijection
        \begin{equation*}
            \operatorname{Aut}(Z,\xi) \overset{1:1}{\longleftrightarrow}
            N_{\pi_{1}(X)}(H)/H\,,
        \end{equation*}
        where $N_{\pi_{1}(X)}(H) = \Set{g \in \pi_{1}(X)\,|\, gHg^{-1} = H}$
        is the \textbf{normalizer subgroup} of $H$.
        It turns out that $H$ is normal, 
        i.e. $N_{\pi_{1}(X)}(H) = \pi_{1}(X)$, if and only if 
        $G := \pi_{1}(X)/H$ acts transitively on the fibres of $\xi$
        and $X$ is homeomorphic to the quotient $Z/G$;
        in this case one calls $\xi$ a \textbf{Galois cover}.
        
        When studying coverings of a normal variety,
        one may always restrict their attention to the Galois case:
        
        \begin{thm}[Galois closure]\label{GC}
            Given a finite algebraic cover $\xi : Z \to X$
            of quasi-projective varieties, there exists a 
            normal quasi-projective variety $\widetilde{Z}$
            and a finite surjective
            morphism $\eta : \widetilde{Z} \to Z$, 
            called the \textbf{Galois closure}, such that
            both $\eta$ and $\xi \circ \eta$ are finite Galois covering maps.
            Moreover, they branch over the same loci.
            \begin{proof}
                See \parencite[Theorem~{3.7}]{GKP2013}.
            \end{proof}
        \end{thm}

        In conclusion, finite covers of a complex normal variety
        are classified by finite index normal subgroups of $\pi_{1}(X)$.
            Since we will deal with \emph{singular} varieties,
        it will be essential to compare (branched) covers of normal varieties with the covers over their smooth loci.
        Indeed, Zariski's Main Theorem allows us to
        extend maps from the smooth locus in a unique way: 

        \begin{thm}[Zariski's Main Theorem]\label{ZMT}
            Any quasi-finite morphism  $\xi : U \to V$ between 
            normal quasi-projective varieties has a unique
            factorization (up to unique isomorpism) into
            an open immersion $\iota : U \hookrightarrow Z$
            into a normal quasi-projective variety $Z$,
            followed by a finite morphism $\widetilde{\xi} : Z \to V$.
            \begin{equation*}
                \begin{tikzcd}
                    U \ar[r, hook, "\iota"] \ar[dr, "\xi"'] & Z \ar[d, "\widetilde{\xi}"] \\
                    & V
                \end{tikzcd}
            \end{equation*}
            Morever, if there exists a group $G$ acting (by algebraic morphisms)
            on $U$ and $V$ such that $\xi$ is $G$-equivariant,
            then the action extends to $Z$ and both $\iota$ and $\widetilde{\xi}$
            are $G$-equivariant.
            \begin{proof}
                See \parencite[Theorem~{3.8}]{GKP2013}.
            \end{proof}
        \end{thm}

        \subsection{The Hilbert scheme of points on a surface}

        Let $\Sigma$ be a normal quasi-projective surface over $\C$;
        normality ensures its singular locus $\textrm{Sing}(\Sigma)$
        has at least codimension $2$ in $\Sigma$.
        The \textbf{Hilbert scheme of $n$ points on $\Sigma$} is the
        set $\operatorname{Hilb}^{n}(\Sigma)$ 
        of $0$-dimensional subschemes $Z \subset \Sigma$
        of length $n$; it turns out this can be given the structure
        of a $2n$-dimensional scheme over $\C$ and has a natural map
        \begin{equation*}
            HC : \operatorname{Hilb}^{n}(\Sigma) \longrightarrow \Sigma^{(n)}\,, \quad
            Z \mapsto \supp Z\,,
        \end{equation*}
        called the \textbf{Hilbert-Chow morphism}.
        We will be interested in the case of $n=2$ points:
        if $\Sigma$ is a smooth surface, it has been proved by \textbf{Fogarty} 
        in \parencite[Theorem~{2.4}]{FOGARTY} that $\hilb{\Sigma}$ is itself a smooth variety;
        in case $\Sigma$ has at worst rational double points singularities, then
        $\operatorname{Hilb}^{2}(\Sigma)$ is an irreducible variety
        with rational Gorenstein singularities 
        by \textbf{Zheng} \parencite[]{ZHENG};
        in these cases, the Hilbert scheme of two points on the surface $\Sigma$ 
        is isomorphic to the following construction:
        call by $\Delta_{\Sigma} \subset \Sigma \times \Sigma$ the diagonal (with the reduced scheme structure) and
        by $\Delta^{(2)}_{\Sigma}$ its image along the quotient 
        $\Ss : \Sigma \times \Sigma \to \Sigma^{(2)} := (\Sigma \times \Sigma)/\xS_{2}$ by the switch of the coordinates;
        the \emph{symmetric product} $\Sigma^{(2)}$ is a complex variety
        such that $\Delta^{(2)}_{\Sigma} \subset \operatorname{Sing}(\Sigma^{(2)})$.
        Let $b$ denote the blow up along $\Delta_{\Sigma}^{(2)}$:
        \begin{equation*}
            \begin{tikzcd}
                         & \Sigma \times \Sigma \arrow[d, "\Ss"] \arrow["\xS_{2}"', loop, distance=2em, in=175, out=125] & \Delta_{\Sigma} \arrow[l, hook] \arrow[d]           \\
\operatorname{Blow}_{\Delta^{(2)}_{\Sigma}}(\Sigma^{(2)}) \arrow[r, "b"] & \Sigma^{(2)}                                                                            & \Delta_{\Sigma}^{(2)} \arrow[l, hook]            \,.
\end{tikzcd}
        \end{equation*}
        Points on $\Sigma^{(2)}$ will be denoted by $p + q := \Ss(p,q) = \Ss(q,p)$ and,
        in particular, %if $p=q$ then 
        $2p := \Ss(p,p)$.

        \begin{lemma}
            Let $\Sigma$ be a normal quasi-projective surface with at worst
            canonical singularities. Then the Hilbert-Chow morphism
            \begin{equation*}
                HC: \operatorname{Hilb}^{2}(\Sigma) \longrightarrow \Sigma^{(2)}
            \end{equation*}
            coincides with the blowing up
            \begin{equation*}
                b: \operatorname{Blow}_{\Delta^{(2)}_{\Sigma}}\left(\Sigma^{(2)}\right) \longrightarrow \Sigma^{(2)}\,.
            \end{equation*}
            \begin{proof}
                As $\Sigma$ is a closed subscheme of $\C^{3}$, then there is a closed immersion
                $\hilb{\Sigma} \hookrightarrow \hilb{\C^{3}}$. Since $\C^{3}$ is smooth,
                then 
                $$\hilb{\C^{3}} \simeq \operatorname{Blow}_{\Delta_{\C^{3}}^{(2)}}(\C^{3})^{(2)}$$
                (for a detailed account, see \parencite[]{CHANG}).
                Then, the inclusion of $\operatorname{Blow}_{\Delta_{\Sigma}^{(2)}}(\Sigma^{(2)})$
                into $\hilb{\C^{3}}$ consists of those $0$-dimensional
                subschemes in $\Sigma$, thus by the universal property of the Hilbert
                scheme, $b$ factors through a morphism $$f:\operatorname{Blow}_{\Delta^{(2)}_{\Sigma}}(\Sigma^{(2)}) \to \hilb{\Sigma}\,,
                \quad b = HC \circ f \,.$$
                Now, since both $b$ and $HC$ induce isomorphisms over 
                $\Sigma^{(2)} \setminus \Delta^{(2)}$, then $f$ is birational.
                As both the domain and range are normal irreducible varieties,
                one concludes $f$ is an isomorphism.
            \end{proof}
        \end{lemma}

        Thus, when talking about the Hilbert scheme $\hilb{\Sigma}$,
        we are going to consider it as the blowing up described above.

    \section{Coverings of the Hilbert scheme}\label{cover}

                The aim of this section is to present a construction
        for quasi-\'etale covers of $\hilb{\Sigma}$,
        induced by quasi-\'etale covers of $\Sigma$,
        as much explicit as possible.
        
        First, let us consider the case $S$ is a smooth surface\footnote{The surface is denoted by $S$ as  `\emph{smooth}',
        to distinguish it from the singular one $\Sigma$.}.
        Then $\hilb{S}$ is smooth as well, and any finite quasi-étale map to it is in fact étale. One might expect to compare its covers 
        to those of $S$ in light of the following:

        \begin{thm}\label{pi-hilb}
            Let $S$ be a smooth (not necessarily compact) complex surface. 
            There exists an isomorphism
            \begin{equation*}
                \pi_{1}(\operatorname{Hilb}^{n}(S)) \simeq
%                \pi_{1}\left(\hilb{S}\right) \simeq 
                \frac{\pi_{1}(S)}{[\pi_{1}(S),\pi_{1}(S)]}\,,
            \end{equation*}
            where the right hand side is the \emph{abelianization} of the fundamental
            group of $S$.
            \begin{proof}
                This fact is due to \textbf{Beauville}, who proved it in
                \parencite[Lemme II.6.1]{BEAUVILLE1983};
                a full detailed proof for $n=2$ can 
                be found in \parencite[Lemma~{7.9}]{GPP-SSS}. 
                Here we rewrite the proof for any $n$ for convenience:
                    consider the big diagonal $\Delta = \bigcup_{1 \le i < j \le n} \Delta_{i,j}$
                    given by 
                    $$\Delta_{i,j} := \Set{ (x_{1}, \dots, x_{n}) \in S^{n} \, | \, x_{i} = x_{j} }\,;$$ 
                    notice $\Delta$ has complex codimension $2$ in $S^{n}$,
                    thus real codimension $4$.
                    Set
                        $$S^{n}_{*} := \Set{ (x_{1}, \dots,x_{n}) \in S^{n} \, | \, \#\{x_{i}\} \ge n-1  }$$
                    the subscheme of $S^{n}$ consisting of those $n$-tuples
                    with \emph{at most} two equal coordinates,
                    thus set $S^{(n)}_{*} := S^{n}_{*}/\xS_{n}$ and $S^{[n]}_{*} := HC^{-1}(S^{(n)}_{*})$. By \parencite[6.(f)]{BEAUVILLE1983} one can
                    identify $S^{[n]}_{*} \simeq \operatorname{Blow}_{\Delta}(S^{n}_{*})/\xS_{n}$.
                    By \parencite[Theorem 2.3, Chapter X]{GODBILLON}, we know that removing subset of high codimension
                    does not affect the fundamental group; 
                    in particular, the inclusion $S^{[n]}_{*} \subset S^{[n]}$
                    induces an isomorphism 
                        $$\pi_{1}(S^{[n]}_{*}) \simeq \pi_{1}(S^{[n]})$$
                    because it has real codimension $4$. 
                    
                    We now study the left hand side by applying \textbf{Seifert-Van Kampen Theorem}: given $U \subset S^{[n]}_{*}$ an open tubular neighborhood
                    of the exceptional divisor $E_{*} \subset S^{[n]}_{*}$
                    and fixed $p \in U \setminus E_{*}$, one has the
                    pushout diagram of groups
                    \begin{equation*}
\begin{tikzcd}
{\pi_{1}(U \setminus E_{*},p)} \arrow[d] \arrow[r] & {\pi_{1}(S_{*}^{[n]} \setminus E_{*},p)} \arrow[d] \\
{\pi_{1}(U,p)} \arrow[r]                           & {\pi_{1}(S^{[n]}_{*},p)}    \,.                      
\end{tikzcd}
                    \end{equation*}
                    By the properties of the blow-up, the Hilbert-Chow morphism induces
                    an isomorphism $S_{*}^{[n]} \setminus E_{*} \simeq S^{(n)}_{*} \setminus \Delta^{(n)} \simeq (S^{n} \setminus \Delta)/\xS_{n}$ with the
                    \emph{unordered configuration space} of $n$ points on $S$;
                    in particular the quotient 
                    $\Ss:S^{n} \setminus \Delta \to (S^{n} \setminus \Delta)/\xS_{n}$, 
                    %since it is a $\xS_{n}$-bundle,
                    for some $p' \in S^{n} \setminus \Delta$ such that
                    $\Ss(p') = HC(p)$,
                    induces the exact sequence of groups
                    \begin{equation}\label{braid}
                        \begin{tikzcd}
                        1 \ar[r] &
                        \pi_{1}(S^{n} \setminus \Delta, p') \ar[r, "\Ss_{*}"] &
                        \pi_{1}(S^{(n)} \setminus \Delta^{(n)},  HC(p)) \ar[r, "\nu"] &
                        \xS_{n} \ar[r] %\ar[l, "s", bend left, dashed]
                        & 1 \,,
                        \end{tikzcd}
                    \end{equation}
                    where $\nu$ sends a class of $\gamma : [0,1] \to (S^{n} \setminus \Delta)/\xS_{n}$ to the permutation $\nu([\gamma])$ such that
                    any lift $\widetilde{\gamma}$ to the ordered configurations $S^{n} \setminus \Delta$ satisfies
                    \begin{equation*}
                        \widetilde{\gamma}_{j}(1) = \widetilde{\gamma}_{\nu([\gamma])(j)}(0)\,,
                        \quad \text{for all } 1 \le j \le n\,.
                    \end{equation*}
                    In fact, the symmetric group acts on 
                    $\pi_{1}(S^{n} \setminus \Delta, p')$
                    by sending the class of a path $\alpha = (\alpha_{1}, \dots,\alpha_{n})$ in $S^{n} \setminus \Delta$ to the class of
                    \begin{equation*}
                        \sigma \cdot \alpha(t) := (\alpha_{\sigma(1)}(t), \dots, \alpha_{\sigma(n)}(t))\,
                    \end{equation*}
                    and moreover the sequence \eqref{braid} splits: 
                    for $1 \le i \le n-1$,
                    let $\eta_{i} : [0,1] \to S$ be a path starting at $p'_{i}$ and ending at $p'_{i+1}$; after setting
                    \begin{equation*}
                        g_{i} : [0,1] \longrightarrow S^{n}\,, \quad g_{i}(t) := (p'_{1},
                        \dots, p'_{i-1}, \eta_{i}(t), p'_{i+1}, \dots, p'_{n})\,,
                    \end{equation*}
                    one may define the path $((i\,,i+1) \cdot g_{i}^{-1}) \ast g_{i}$
                    which exchanges the $i$-th and the $(i+1)$-th coordinates of the base point $p'$; by dimension reasons, one may consider a homotopic path $\tau_{i}$ which avoids $\Delta_{i,i+1}$, i.e. there exists
                    $\tau_{i} : [0,1] \to S^{n} \setminus \Delta$ and
                    $\tau_{i} \simeq ((i\,, i+1) \cdot g_{i}^{-1}) \ast g_{i}$
                    in $S^{n}$.
                    Notice $\tau_{i}$ is not a loop in $S^{n} \setminus \Delta$,
                    but its image $\Ss(\tau_{i})$ is a loop ``winding around'' $\Ss(\Delta_{i,i+1})$ and thus one can define  
                    a section $s$ on transpositions by 
                        \begin{equation*}
                            s(i,i+1) := [\Ss(\tau_{i})]\,, \quad \text{for all } 1 \le i \le n-1\,,
                        \end{equation*}
                    hence the $\xS_{n}$ acts on $[\lambda] \in \pi_{1}(S^{n} \setminus \Delta, p')$ as
                    \begin{equation*}
                        (i,i+1) \cdot [\lambda] := [\tau_{i}^{-1} \ast ((i,i+1) \cdot \lambda) \ast \tau_{i}]
                    \end{equation*}
                     Now, since $\pi_{1}(S^{n} \setminus \Delta) \simeq \pi_{1}(S^{n})$, by putting everything together one deduces
                    \begin{align*}
                        \pi_{1}(S_{*}^{[n]} \setminus E_{*},p) 
                        &\simeq \pi_{1}(S^{n},p') \rtimes \xS_{n} \\
                        &\simeq \pi_{1}(S^{n},(p_{0}, \dots,p_{0})) \rtimes \xS_{n}
                        % (\pi_{1}(S,p'_{1}) \times \dots \times \pi_{1}(S,p'_{n})) \rtimes \xS_{n}
                        \simeq \pi_{1}(S,p_{0})^{n} \rtimes \xS_{n}\,,
                    \end{align*}
                    where the isomorphism in the middle is induced
                    by a path from $p'$ to some  $(p_{0}, \dots, p_{0})$ in $S^{n}$.
                    In this case the action of $\xS_{n}$ translates on $\pi_{1}(S,p_{0})^{n}$
                    as the permutation of the entries.

                    In the previous construction, 
                    one may take a representative $\tau_{i}$ ``close enough''
                    to the diagonal $\Delta_{i,i+1}$ in such a way that its
                    image is contained in a tubular neighborhood mapped into $U$,
                    so that $[\Ss(\tau_{i})]$ can be seen as an element in $\pi_{1}(U \setminus E_{*}, p)$; the inclusion $U \setminus E_{*} \subset U$ sends
                    all these $[\Ss(\tau_{i})]$ to the identity in $\pi_{1}(U,p)$,
                    and hence by Seifert-Van Kampen the group $K$ generated by these loops
                    imposes relations in $\pi_{1}(S^{[n]}_{*},p)$; by taking track
                    of all the isomorphisms above, one may identify
                    \begin{equation*}
                        K \simeq 
                        \big\langle (\sigma \cdot g^{-1}) \ast g\, | \, \sigma \in \xS_{n}, g \in \pi_{1}(S,p_{0})^{n} \big\rangle
                    \end{equation*}
                    and thus we obtain the short exact sequence
                    \begin{equation*}
                        \begin{tikzcd}
                        1 \ar[r] &
                        \langle (\sigma \cdot g^{-1}) \ast g \rangle \ar[r, "j"] &
                        \pi_{1}(S,p_{0})^{n} \rtimes \xS_{n} \ar[r] &
                        \pi_{1}(S^{[n]}_{*},p) \ar[r] & 1 \,
                        \end{tikzcd}
                    \end{equation*}
                    as noted by Beauville.

                    It remains to prove that the cokernel of $j$ is isomorphic
                    to the abelianization of $\pi_{1}(S,p_{0})$:
                    the map which composes a $n$-tuple of paths
                    \begin{align*}
                        \ast : \pi_{1}(S,p_{0})^{n} \rtimes \xS_{n} \longrightarrow 
                        \frac{\pi_{1}(S,p_{0})}{[\pi_{1}(S,p_{0}), \pi_{1}(S,p_{0})]}\,,
                        \\ \big(g_{1}, \dots, g_{n}, \sigma)
                        \mapsto [g_{\sigma(1)} \ast g_{\sigma(2)} \ast \dots \ast g_{\sigma(n)}]
                    \end{align*}
                    is well defined because the target is abelian, moreover it is surjective.
                    Notice that the relations in $K$ imply that for any 
                    $[g] = [(g_{1}, g_{2}, \dots, g_{n})]$ in the cokernel of $j$ one has the identity $[\sigma \cdot g] = [g]$ for any permutation $\sigma \in \xS_{n}$, thus
                    \begin{align*}
                        [(g_{1}, g_{2}, \dots, g_{n})] 
                        &= [(g_{1}, 1, g_{3}, \dots, g_{n})] \ast [(1, g_{2}, 1,\dots, 1)] \\ 
                        &= [(g_{1}, 1, g_{3}, \dots, g_{n})] \ast [(g_{2}, 1,\dots, 1)] \\
                        &= [(g_{1} \ast g_{2}, 1, g_{3}, \dots, g_{n})] \\
                        &= [(g_{1} \ast g_{2} \ast g_{3}, 1, 1, g_{4}, \dots, g_{n})] \\
                        &= \dots = [(\ast(g), 1, 1, \dots, 1)] 
                    \end{align*}
                    from which we deduce the above morphism is injective, hence the desired isomorphism.
            \end{proof}
        \end{thm}

    \subsection{Construction of the coverings}\label{AUT}
        A consequence of \textbf{Theorem~\ref{pi-hilb}} is that $\pi_{1}(\hilb{S})$
        has, at most, as many (isomorphism classes of) Galois \'etale coverings %Galois (\textcolor{red}{Anche tu lo noti quindi diamo una bella decisione. Per questo teorema tutti i sottogruppi di $\pi_1(\hilb{S})$ sono normali in quanto abeliano. Quindi TUTTI i rivestimenti di $\hilb{S}$ sono Galois. Quello che direi io è che visto che uno è un quoziente dell'altro per ogni rivestimento (Galois o meno) di $\hilb{S}$ esiste (almeno) un rivestimento di $S$ con lo stesso indice}) coverings
        as the surface $S$, so it is natural to ask if there exists a geometric description 
        of the way coverings are induced to the Hilbert scheme.
        We now present an explicit construction for these covers based on the one made for Enriques surfaces by \textbf{Taro Hayashi} in \cite{TARO}. 
        
        Let $\xi:Z \to S$ be a finite  Galois \'etale covering
        with \emph{abelian} deck transformation group $G$;
        assume its order is $\lvert G \rvert = d$.
        Notice that the product 
        $$\xi^{2}:Z \times Z \to S \times S$$
        is again a covering and the symmetric group $\xS_{2}$
        generated by the transposition $\sigma = (1 \,, 2)$
        acts on $\xi^{2}$, i.e. $\xi^{2}$ is a $\xS_{2}$-equivariant
        morphism, thus the swap of coordinates is an automorphism of the covering.
        Also the group $G$ induces automorphisms of $\xi^{2}$:
        for $g \in G$ and $(z,w) \in Z^{2}$, set
        \begin{equation*}
            s_{g}(z,w) := (z, g \cdot w)\,, \quad 
            \delta_{g}(z,w) := (g \cdot z, g \cdot w)\,.
        \end{equation*}
        Let now $J \le \operatorname{Aut(Z^{2})}$ 
        (resp. $H \le \operatorname{Aut(Z^{2})}$)
        be the subgroup generated by 
        $\xS_{2} = \langle \sigma \rangle$ and $\Set{s_{g}|g \in G}$
        (resp.
        by $\xS_{2}$ and $\Set{\delta_{g}|g \in G}$). 
        After checking that, for any $g,g',h\in G$ there hold the relations
        \begin{equation*}
            s_{h}\sigma s_{g'}\sigma s_{g}=\sigma s_{g'} \sigma s_{hg}\,, \quad
            s_{g}\sigma s_{g'} \sigma = \sigma s_{g'} \sigma s_{g}\,, 
        \end{equation*}
        one notices that any element $j \in J$ can be written in the form
        $j=\tau's_{g'} \tau s_{g}$, for some $g,g'\in G$ and $\tau, \tau'\in \xS_{2}$; moreover, there is a well-defined sign homomorphism
        \begin{equation*}
            \operatorname{sgn} : J \to \xS_{2}\,, \quad
            \operatorname{sgn}(\tau's_{g'} \tau s_{g}) := \tau'\tau\,,
        \end{equation*}
        which sits in a split sequence
        \begin{equation}\label{GxG}
                        \begin{tikzcd}
                        0 \ar[r] &
                        G \times G \ar[r, "\xG"] &
                        J \ar[r, "\operatorname{sgn}"] &
                        \xS_{2} \ar[r] \ar[l, "s", bend left, dashed]
                        & 1 \,,
                        \end{tikzcd}
        \end{equation}
        where $\xG(g',g) := \sigma s_{g'} \sigma s_{g}$ and $s(\tau) = \tau$.
        From the equation
        \begin{equation*}
            \delta_{g} = \sigma s_{g} \sigma s_{g}\,, \quad \text{for all } g \in G\,,
        \end{equation*}
        we deduce that $H$ is contained in $J$ and hence $H \le J$; 
        moreover, $H$ is normal in $J$ because $G$ is abelian:
        indeed, it is clear that $\sigma \delta_{g}\sigma = \delta_{g}$ for any $g \in G$, while for $h \in G$ one has:
        \begin{equation*}
            s_{h}^{-1}\delta_{g}s_{h}(z,w) = (g \cdot z, h^{-1}gh \cdot w) = 
            (g \cdot z, gh^{-1}h \cdot w) = \delta_{g}(z,w).
        \end{equation*}
        
        From the sequence \eqref{GxG}, one deduces $Z^{2}/J \simeq S^{(2)}$
        and hence there is a commutative triangle
        \begin{equation}\label{q}
            \begin{tikzcd}
                                         & Z^{2} \arrow[dd, "\Xi"] \arrow[ld, "q"'] \\
Z^{2}/H \arrow[rd, "\widetilde{\xi}"'] &                                      \\
                                         & S^{(2)}        \,.                     
\end{tikzcd}
        \end{equation}
        One may describe the fibres of $\Xi$  set-theoretically  as
        \begin{equation*}
            \Xi^{-1}(\xi(z) + \xi (w)) = 
            \Set{(g \cdot z, h \cdot w), (h \cdot w, g \cdot z) \, | \, g,h \in G }\,, 
        \end{equation*}
        which has cardinality $2d^{2}$ whenever $\xi(z) \ne \xi(w)$,
        while the cardinality drops to $d^{2}$ when $\xi(z) = \xi(w)$.

        For each $g \in G$, set 
        \begin{equation*}
            T_{g} := \Set{(z,g \cdot z) \, | \, z \in Z, g \in G} \subset Z^{2}\,;
        \end{equation*}
        then $\Xi^{-1}(\Delta^{(2)}) = \cup_{g \in G} \, T_{g} =: T$
        is a closed $J$-invariant subscheme of codimension $2$ and, 
        after removing it, the resulting restriction
        \begin{equation}\label{outside-exceptional}
            \Xi' : Z^{2} \setminus T \longrightarrow S^{(2)} \setminus \Delta^{(2)}_{S}
        \end{equation}
        is an \'etale covering map of degree $2d^{2}$;
        this shows $\Xi$ is a finite quasi-\'etale cover, 
        which ramifies on $\Delta_{S}^{(2)}$.
        Notice that, for $g,h \in G$, one has $(z,w) \in T_{g} \cap T_{h}$
        if and only if $g \cdot w = h \cdot w$; this means that $h^{-1}g$
        has a fixed point, thus $g=h$ because the $G$ acts freely on $Z$.
        In particular, $\Xi^{-1}(\Delta^{(2)})$ consists of $\lvert G \rvert$
        disjoint copies of the diagonal $\Delta_{Z} = T_{\operatorname{id}}$.

        \begin{lemma}
            If $t \in H$ fixes a point $(z,g \cdot z) \in T_{g}$,
            then it fixes the whole $T_{g}$ pointwise.
            \begin{proof}
                Since for every $g,h \in G$ and $\tau \in \xS_{2}$ one has
                \begin{equation*}
                    \tau \circ \delta_{g} \circ \tau = \delta_{g}\,, \quad
                    \delta_{h} \circ \tau \circ \delta_{g} = \tau \circ \delta_{hg}\,,
                \end{equation*}
                any element in $t \in H$ can by written as 
                $t=\tau \circ \delta_{g}$, 
                for some $\tau \in \xS_{2}\,, h \in G$.
                Thus there are two cases:
                \begin{itemize}
                    \item if $\tau=\operatorname{id}$, 
                    then $(z, g \cdot z) = (h \cdot z, hg \cdot z)$,
                    which implies $h=\operatorname{id}$;
                    \item if $\tau = \sigma$, 
                    then $(z, g \cdot z) = (hg \cdot z,h \cdot z)$
                    implies $hg=\operatorname{id}$ and $g=h$,                   thus $g^{2}=\operatorname{id}$. In this case,
                    for every $z \in Z$ one has
                    \begin{equation*}
                        t(z,g\cdot z) = (g^{2} \cdot z, g \cdot z) 
                        = (z, g \cdot z)\,.
                    \end{equation*}
                \end{itemize}
                In particular, only those components $T_{g}$ 
                corresponding to $\operatorname{ord}(g)\le2$ remain fixed.
            \end{proof}
        \end{lemma}

        In particular, this means that, 
        when restricted to
        \begin{equation*}
            \Xi\vert_{T} : T \longrightarrow \Delta_{S}^{(2)}\,,
        \end{equation*}
        the map $\Xi$ induces a covering of the ramification locus.
            
        \begin{lemma}
            The cardinality of the fibres of the quotient 
            $\widetilde{\xi}:Z^{2}/H \to S^{(2)}$ 
            is constant, equals $\lvert G \rvert$.
            \begin{proof}
                Given $(z,w) \in Z^{2}$, 
                then its $H$-orbit consists of 
                \begin{equation*}
                    H \cdot (z,w) = \Set{(g \cdot z, g \cdot w), (g \cdot w, g \cdot z) \, | \, g \in G }\,,
                \end{equation*}
                thus any point in the fibre 
                $\widetilde{\xi}^{-1}(\xi(z) + \xi (w))$
                admits a representative of the form $(z, g \cdot w)$,
                for some $g \in G$. Now, if $g,h \in G$, then the equation
                \begin{equation*}
                    q(z,g \cdot w) = q(z,h \cdot w)
                \end{equation*}
                implies there exists $k \in G$ such that
                either $(z, g \cdot w) = (k \cdot z, kh \cdot w)$
                or $(z, g \cdot w) = (kh \cdot w, k \cdot z)$.
                In the former case, $k \cdot z = z$ implies $k = \operatorname{id}_{Z}$ because the action of $G$ is free on $Z$,
                thus one deduces also $w = g^{-1}h \cdot w$, so $g=h$.
                The latter case implies 
                $(z,g \cdot w) = (z, k \cdot z) \in T_{k}$
                and hence we see $q(z,g \cdot z) = q(z,k \cdot z)$
                if and only if $(z, g \cdot z) \in T_{k}$, which implies $g=k$.
                We have thus proven that 
                $\Set{(z,g \cdot w)| g \in G}$
                consists of non-equivalent representatives.
            \end{proof}
        \end{lemma}

        As a consequence, we obtain that the cover
        \begin{equation*}
            \widetilde{\xi} : Z^{2}/H \longrightarrow S^{(2)}
        \end{equation*}
        is unramified. Since the blow-up along $\Delta_{S}$ 
        is a resolution of the singularities of $S^{(2)}$
        and $q(T) = \widetilde{\xi}^{-1}(\Delta_{S}^{(2)})$,
        by the universal property of the blow-up we obtain 
        a smooth and unramified map $\xi^{[2]}$ of degree 
        $\lvert G \rvert$:
        \begin{equation}\label{xi-hilb}
            \begin{tikzcd}
            \operatorname{Blow}_{q(T)}\left(\frac{Z^{2}}{H}\right)  \arrow[d, "\exists !"', "\xi^{[2]}"] \arrow[r, "t"] 
            & Z^{2}/H \arrow[d, "\widetilde{\xi}"]  & Z^{2}\arrow[l] \arrow[dl, "\Xi"] \\
            \operatorname{Blow}_{\Delta^{(2)}}(S^{(2)})  \arrow[r, "HC"]                         & S^{(2)}    &           \,.     
            \end{tikzcd}
        \end{equation}

        In conclusion, the assignment $\xi \mapsto \xi^{[2]}$ is the geometric realisation we were aiming for and we can prove the following theorem.
        \begin{thm}\label{teorema smooth}
            Let $S$ be a smooth complex surface with finite fundamental group
            and $\eta:\xY \to \hilb{S}$ be a finite \'etale covering. 
            Then, using the notation of \Cref{xi-hilb}, 
            there exists an étale covering $\xi:Z \to S$ such that
            $\eta$ is isomorphic to the covering $\xi^{[2]}$.
            \begin{proof}
                Recall that $\eta$ corresponds to a subgroup 
                $N \le \pi_{1}(\hilb{S})$ and, since the latter is abelian,
                then $N$ is normal and hence the covering $\eta$ is Galois,
                with abelian deck transformations group $G := \pi_{1}(\hilb{S})/N$.
                Moreover, by \textbf{Theorem~\ref{pi-hilb}} we know
                the fundamental group of $\hilb{S}$ is the abelianization
                of $\pi_{1}(S)$, hence there exists a normal subgroup $N'\lhd \pi_{1}(S)$ such that $G \simeq \pi_{1}(S)/N'$; this
                corresponds to a finite \'etale cover $\xi:Z \to S$
                and hence, by the explicit construction described in this Section, this induces
                the cover $\xi^{[2]}$ on $\hilb{S}$, 
                which has still deck transformation group $G$,
                implying that $\eta$ is isomorphic to $\xi^{[2]}$ as a covering. 
            \end{proof}
        \end{thm}

\subsection{Approximating \'etaleness}

    Let now $\Sigma$ be a quasi-projective surface over $\C$
    with, at worst, ADE singularities and denote
    by $S := \Sigma^{sm}$ its smooth locus. Normality ensures
    the singular locus $\Sigma \setminus S$ has complex codimension $2$ in $\Sigma$,
    thus it consists of finitely-many isolated points.
    In this subsection we would like to \emph{characterize} quasi-\'etale maps
    \begin{equation*}
        \pi:\Zz \longrightarrow \hilb{\Sigma}\,.
    \end{equation*}
    By describing the Hilbert scheme as a \emph{quiver variety},
    \textbf{Craw} and \textbf{Yamagishi} proved in \parencite[Corollary~{6.6}]{CRAW-YAMAGISHI} that $\hilb{\Sigma}$ is a normal variety. 
    Thus, by purity we know that $\pi$ branches along the singular locus of $\hilb{\Sigma}$,
    which still has codimension $2$; thus, $\pi$ is quasi-\'etale
    if and only if
    \begin{equation*}
        \pi^{sm} : \Zz^{\circ} \longrightarrow \hilb{\Sigma}^{sm}\,,
    \end{equation*}
    is \'etale, where $\Zz^{\circ} := \pi^{-1}(\hilb{\Sigma}^{sm})$.
    Now, one might wonder if this relates to the coverings of $\hilb{S}$
    somehow, for $\hilb{S}$ can be embedded in $\hilb{\Sigma}^{sm}$;
    notice this embedding is neither open, nor closed in general,
    for one may describe
    \begin{equation*}
        \hilb{\Sigma}^{sm} \setminus \hilb{S}
    \end{equation*}
    rather explicitly, as done in the proof
    of \parencite[Proposition~{7.10}]{GPP-SSS}: 
    let $\textrm{Sing}(\Sigma) = \{p_{1}, \dots, p_{k}\}$
    and $\Ss: \Sigma \times \Sigma \to \Sigma^{(2)}$ be the symmetric quotient; 
    if we denote by $Z_{i}:=\Ss(\{p_{i}\} \times \Sigma) = \Ss(\Sigma \times \{p_{i}\})$
    and by $2p_{i} :=  \Ss(p_{i},p_{i})$;
    then $$\textrm{Sing}\left(\Sigma^{(2)}\right) = \Ss(\Delta_{\Sigma}) \cup Z_{1} \cup \dots \cup Z_{k}$$ 
    and, by putting $W_{i} := HC^{-1}(Z_{i} \setminus \{2p_{i}\})$,
    one deduces that
    \begin{equation*}
        \textrm{Sing}\left(\hilb{\Sigma}\right) = \overline{W_{1}} \cup \dots \cup \overline{W_{k}}\,.
    \end{equation*}
    Notice that for each component 
    $\overline{W_{i}} \setminus HC^{-1}(2p_{i}) \simeq Z_{i} \setminus \{2p_{i}\}$
    is isomorphic to $\Sigma \setminus \{p_{i}\}$, thus has codimension $2$ in $\hilb{\Sigma}$;
    it remains to study the fiber of each $2p_{i}$: since each $p_{i}$ is an ADE
    singularity, then the tangent space $T_{p_{i}}\Sigma$ 
    has dimension $3$ for each $i=1, \dots, k$,
    and hence by $HC^{-1}(2p_{i}) \simeq \P(T_{p_{i}}\Sigma)$ we conclude that 
    $HC^{-1}(2p_{i})$ is isomorphic to a $\P^{2}$. Then $\hilb{\Sigma}^{sm}$
    differs from $\hilb{S}$ by the complements $\P(T_{p_{i}}\Sigma) \setminus \overline{W_{i}}$,
    which are $2$-dimensional.

    After setting $\Zz' := \pi^{-1}(\hilb{S})$, define
    \begin{equation*}
        \pi':\Zz'\longrightarrow \hilb{S}\,.
    \end{equation*}

    \begin{lemma}\label{apx}
        Given $\pi:\Zz \longrightarrow \hilb{\Sigma}$, then
        $\pi^{sm}$ is \'etale if and only if $\pi'$ is \'etale.
        \begin{proof}
            Since \'etaleness is stable under base change,
            then it is clear that $\pi^{sm}$ \'etale implies
            that its restriction $\pi'$ is \'etale too. Conversely,
            if $\pi^{sm}$ is not \'etale, then by purity
            it must ramify along a divisor; since 
            $\hilb{\Sigma}^{sm} \setminus \hilb{S}$ has at most
            dimension $2$, then the branch locus must intersect
            $\hilb{S}$, thus one concludes $\pi'$ is not \'etale.
        \end{proof}
    \end{lemma}

    Thanks to the results obtained in the smooth case,
    it seems more manageable to build covers of the smooth
    locus $\hilb{\Sigma}^{sm}$ and hope to extend these to  
    (possibly branched) covers to the whole Hilbert scheme.
    \begin{cor}
        Let $V$ be a normal, complex, quasi-projective variety. 
        Then quasi-\'etale finite covers of  $V$ 
        are in one-to-one correspondence with finite \'etale covers of $V^{sm}$.
        \begin{proof}
            Since all varieties involved are normal,
            then by purity of the branch locus it follows that any
            finite quasi-\'etale cover $f:\Zz \to V$ branches over 
            $\textrm{Sing}(V)$, thus $f$ restricts to the \'etale cover
            $f\vert_{f^{-1}(V^{sm})}:f^{-1}(V^{sm}) \to V^{sm}$.
            Conversely, if $g : \Yy \to V^{sm}$ is finite and \'etale, 
            then by Theorem~\ref{ZMT} the composite $\Yy \overset{g}{\to} V^{sm} \subset V$
            factorizes through a (unique) finite cover $\widetilde{\Yy} \to V$.
        \end{proof}
    \end{cor}

    We are now ready to prove the analogue of \textbf{Theorem~\ref{teorema smooth}}
    in the singular setting; in the statement, we use the same notations
    as in the previous subsection and adapt it for singular surfaces too.
        \begin{thm}\label{sing-hilb-cover}
            Let $\Sigma$ be a normal, complex, quasi-projective surface with ADE singularities at worst. 
            Then, given a finite quasi-\'etale covering $\xi:W \to \Sigma$,
            the map $\xi^{[2]}$ in \eqref{xi-hilb}
            is a quasi-\'etale covering of $\hilb{\Sigma}$. Moreover, all finite quasi-\'etale coverings of $\Sigma$ arise in this way.
            \begin{proof}
            
                Let $\xi:W \to \Sigma$ be a finite-quasi \'etale covering,
                which restricts over $S = \Sigma^{sm}$ 
                to a Galois covering $\xi^{sm}$ 
                with group $G$; recall that,  by \textbf{Theorem~\ref{ZMT}}, $G$ acts on $\xi$ as well. 
                If we denote by $H$ the group of automorphisms
                introduced at the beginning of \textbf{Section~\ref{AUT}}
                and consider the quotient map $q:W^{2} \to W^{2}/H$,
                then by purity of the branch locus we know that 
                $$\xi^{[2]}:\operatorname{Blow}_{q(T)}\left(\frac{W^{2}}{H}\right)
                \to \hilb{\Sigma}$$ is quasi-\'etale 
                if and only if $(\xi^{[2]})^{sm}$ is \'etale.
                By \textbf{Lemma~\ref{apx}}, this is equivalent
                to the \'etaleness of the restriction
                \begin{equation}\label{restricted}
                    \pi': \operatorname{Blow}_{T'}\left(\frac{\xi^{-1}(S) \times \xi^{-1}(S)}{H}\right) 
                    \to \hilb{S}\,,
                \end{equation}
                where $T' = \cup_{g \in G} \{q(z, g \cdot z) \,|\, z \in \xi^{-1}(S)\}$
                is the same as $q(T)$, with finitely many punctures. 
                Since the construction is the one described in \textbf{Section~\ref{AUT}} for the smooth case, 
                the map \eqref{restricted} is
                \'etale, hence the first claim follows. 
                The second part follows from \textbf{\Cref{teorema smooth}} applied to the smooth part and
                the fact that a birational morphism between normal varieties
                is, in fact, an isomorphism.
            \end{proof}
        \end{thm}

    Now, it follows directly from \Cref{teorema smooth} and \Cref{sing-hilb-cover} the first theorem we wanted to prove, i.e. \textbf{\Cref{ThmA}}.

    \section{Application to irreducible symplectic surfaces}\label{sec:ISV}
        
In this final section of the paper, 
we will apply the machinery of \textbf{Section~\ref{cover}}
in case $\Sigma$ is any \emph{projective irreducible symplectic surface}
in order to prove \textbf{Theorem~A} 
(that is \textbf{\Cref{thmA}} in the second paragraph).
Before that, we recall in the first subsection the definition of ISV.

\subsection{Irreducible symplectic varieties}

    We now recall the notion of \emph{irreducible symplectic variety},
    a singular analogue of a \emph{hyperk\"ahler manifold}
%    \emph{irreducible symplectic manifolds} 
    proposed by Greb-Kebekus-Peternell in \parencite[Definition~{8.6.2}]{GKP2016}:
    while in the smooth setting these manifolds carry a 
    nowhere-degenerate alternating $2$-form, in the singular setting
    we require the existence of such a structure only
    on the smooth locus, in such a way that it extends
    via the above mentioned quasi-\'etale coverings.
    Let us recall here the definition: given $X$ a complex normal K\"ahler variety,
    let  the \textbf{sheaf of reflexive $p$-forms} on $X$ be
    \begin{equation*}
         \Omega_{X}^{[p]} :=  (\Omega_{X^{sm}}^{p})^{\vee \vee}
         \simeq \iota_{*}  \Omega_{X^{sm}}^{p}\,,
    \end{equation*}
    where $\iota_{*}$ denotes the pushforward along the embedding
    $\iota : X^{sm} \hookrightarrow X$. We can then define
    the \textbf{reflexive Dolbeaut cohomology} of $X$ as
    \begin{equation*}
        H^{[p],q}(X) := H^{q}\left(X, \Omega_{X}^{[p]}\right) \,.
    \end{equation*}

\begin{comment}
    In case $X$ is projective with, at worst, canonical singularities,
    then it holds the \textbf{Hodge duality}
    \begin{equation*}
        H^{0}\left(X, \Omega_{X}^{[p]}\right) 
        \simeq H^{0}\left(X^{sm}, \Omega_{X^{sm}}^{p}\right) 
        \simeq H^{p}(X,\Oo_{X})
    \end{equation*}
    as stated in \parencite[Proposition~{6.9}]{GKP2016};
    if moreover $\omega_{X} \simeq \Oo_{X}$, then
    by \parencite[Proposition~{6.1}]{GKP2016}
    the wedge product gives a non-degenerate pairing
    \begin{equation*}
        \wedge :  H^{[p],0}(X) \times  H^{[n-p],0}(X) \to H^{0}(X, \omega_{X}) \simeq \C\,.
    \end{equation*}
\end{comment}    
    
    \begin{df!}\label{ISV}
        An \textbf{irreducible symplectic variety} (shortly, \textbf{ISV}) is a normal compact 
        K\"ahler space $X$ that verifies the three following conditions:
        \begin{enumerate}
            \item it has canonical singularities,
            \item it admits a holomorphic symplectic form 
            $\sigma$ on its smooth locus $X^{sm}$, that we will view as a global section of
            the sheaf of reflexive $2$-forms $\Omega^{[2]}_{X} = \left(\Omega_{X}^{2}\right)^{\vee \vee}$;
            \item for every \emph{finite quasi-\'etale covering} $f : Y \to X$
            we have that
                \begin{equation*}
                    \bigoplus_{p} H^{0}(Y, \Omega_{Y}^{[p]}) \simeq \C\left[f^{[*]}\sigma\right]\,,
                \end{equation*}
            as a $\C$-algebra, where 
            $f^{[*]}\sigma := \left(f^{*}\sigma \right)^{\vee \vee}$.
        \end{enumerate}
        In case $X$ is $2$-dimensional, then it will be called \textbf{irreducible symplectic surface}.
    \end{df!}
    
    Note that condition $3$ in the above Definition is a 
    mimicking of the external edge of the Hodge diamond
    of a hyperk\"alher manifold, realized via reflexive forms.
    It is known by \parencite[Proposition~{2.19}]{PEREGO-EX} that 
    any ISV $X$ has $h^{1}(X,\Oo_{X})=0$ and $h^{[2],0}(X)=1$.

    \begin{rmk!}
        A finite quasi-\'etale covering $f:Y \to X$ of
        an ISV $X$ is itself an ISV: for quasi-\'etale maps
        preserve canonical singularities and the smooth
        locus of $Y$ is endowed with the form $f^{[*]}\sigma$.
        Then condition 3. holds by naturality of the
        reflexive pullback.
    \end{rmk!}

    \subsection{Covers of symplectic varieties}

        Let now $\Sigma$ be an irreducible symplectic surface.
        In this case, $\operatorname{Sing}(\Sigma) = \Sigma \setminus S$ 
        consists of finitely many singular points of type ADE (see \parencite[Theorem~6]{NAMIKAWA}) and by 
        Xiao's classification \parencite[Theorem~3]{XIAO} we know that
        $\pi_{1}(S)$ is always a finite group, 
        hence we conclude that
        the bijection described in \textbf{Section~\ref{cover}} between
        all finite quasi-\'etale coverings of $\hilb{\Sigma}$ and subgroups of
        \begin{align*}
%            \pi_{1}^{\text{\'et}}\left( \hilb{S} \right)
            \pi_{1}\left( \hilb{S} \right) %^{\wedge} 
%            \simeq \left( \frac{\pi_{1}(S)}{[\pi_{1}(S),\pi_{1}(S)]} \right) ^{\wedge}
            \simeq \frac{\pi_{1}(S)}{[\pi_{1}(S),\pi_{1}(S)]} \,,
        \end{align*}
        holds true. 

        \begin{thm}\label{thmA}
            If $\Sigma$ is an irreducible symplectic surface,
            then $\hilb{\Sigma}$ is an irreducible symplectic variety
            of dimension $4$.
            \begin{proof}
                It is already known that $\hilb{\Sigma}$ is a 
                \emph{primitive symplectic variety}, 
                i.e. $\hilb{\Sigma}$ fulfills conditions $1.$ and $2.$
                in \textbf{Definition~\ref{ISV}},
                for it is proven in \parencite[Corollary~{7.7}]{GPP-SSS};
                thus it is sufficient to prove condition 3. 
                
                Given a finite quasi-\'etale cover $\eta:\xY \to \hilb{\Sigma}$,
                we know by Zariski's Main Theorem that it is determined by $\eta^{sm}$,
                which is Galois by abelianess of $\pi_{1}(\hilb{\Sigma^{sm}})$;
                thus, $\eta$ is Galois with deck transformation group $G$. By \textbf{Theorem~\ref{sing-hilb-cover}} it 
                is induced by a quasi-\'etale Galois  covering 
                $\xi:Y \to \Sigma$, 
                i.e.
                $\eta = \xi^{[2]}$.
                As $\xY$ is the blow-up $t$ in the diagram \eqref{xi-hilb} with $Z=Y$ and $S = \Sigma$, then 
                it is birational to $Y^{2}/H$, where $H$ is the group
                defined in \textbf{Section~\ref{AUT}}; moreover,
                notice $Y^{2}/H \simeq Y^{(2)}/G$, 
                where the action of $G$ on the symmetric product is the diagonal one, i.e.
                $g \cdot (z+w) := g \cdot z + g \cdot w$.
%                $\Sigma \times \Sigma/H$, where we recall that $H$ is the subgroup of $\operatorname{Aut}(Y \times Y)$ generated by $\xS_{2}$ and $G$. 
                Since by \parencite[Lemma~{7.5}]{GPP-SSS} 
                the external reflexive Hodge numbers are birational invariants, 
                for every $p \in \Z$ it holds:
                \begin{equation*}
                    h^{[p],0}(\xY) % \simeq h^{[p],0}(Y \times Y/H) 
                    = h^{[p],0}(Y^{(2)}/G) = h^{[p],0}\left(Y^{(2)}\right)^{G}\,;
%                    \simeq h^{[p],0}(Y \times Y)^{H}\,;
                \end{equation*}
%                where we used that $H \simeq G \rtimes \xS_{2}$\todo{Qua ho fatto un po' di parkour, non so se sia vero...};
                on the other hand, $Y^{(2)}$ inherits a reflexive $2$-form form from $\Sigma^{(2)}$ (for it is a PSV) and $G$ acts symplectically on it, implying
                that
                \begin{equation*}
                    h^{[p],0}\left(Y^{(2)}\right)^{G} = h^{[p],0}\left(\Sigma^{(2)}\right) 
                    = h^{[p],0}\left(\hilb{\Sigma}\right)\,.
                \end{equation*}
                By putting everything together, we conclude that 
                $h^{[p],0}(\xY) = h^{[p],0}(\hilb{\Sigma})$
                for every $p \in \Z$; in particular, $H^{[2],0}(\xY)$
                is $1$-dimensional and must be generated by the reflexive pullback $f^{[\ast]}\sigma$ of the generator $\sigma \in H^{2,0}(\hilb{S})$.
            \end{proof}
        \end{thm}

        \begin{rmk}
            Notice that by the existence of a (reflexive) symplectic form on $\hilb{\Sigma}$, we deduced \emph{a posteriori} that $G$ acts symplectically on its covering $\xY$.
            In case $\Sigma$ did not have any $2$-form
            and some covering space $W$ was symplectic,
            then we would deduce that $G$ does not preserve any symplectic form;
            in particular, by construction also the quotient $W^{2}/H$
            would not have any symplectic $2$-form
            and hence $\hilb{\Sigma}$ is a variety of Calabi-Yau type, 
            as in the case of \textbf{Hayashi} \parencite[]{TARO}.
        \end{rmk}

\clearpage
\phantomsection
\addcontentsline{toc}{section}{References}\include{ref.bib}

\printbibliography[title={References}]

\end{document}